Computer-based and paper-and-pencil tests: A study in calculus for STEM majors

Lawrence Smolinsky[1], Brian D. Marx[2], Gestur Olafsson[1], and Yanxia A. Ma[2]

[1] Depart of Mathematics, Louisiana State University, Baton Rouge, LA 70803 USA
smolinsk@math.lsu.edu

[2] Depart of Experimental Statistics, Louisiana State University, Baton Rouge, LA 70803

## Abstract

Computer-based testing is an expanding use of technology offering advantages to teachers and students. We studied Calculus II classes for STEM majors using different testing modes. Three sections with 324 students employed: Paper-and-pencil testing, computer-based testing, and both. Computer tests gave immediate feedback, allowed multiple submissions, and pooling. Paper-and-pencil tests required work and explanation allowing inspection of high cognitive demand tasks. Each test mode used the strength of its method. Students were given the same lecture by the same instructor on the same day and the same homework assignments and due dates. The design is quasi-experimental, but students were not aware of the testing mode at registration. Two basic questions examined were: (1) Do paper-and-pencil and computer-based tests measure knowledge and skill in STEM Calculus II in a consistent manner? (2) How does the knowledge and skill gained by students in a fully computer-based Calculus II class compare to students in a class requiring pencil-and-paper tests and hence some paper-and-pencil work. These results indicate that computer-based tests are as consistent with paper-and-pencil tests as computer-based tests are with themselves. Results are also consistent with classes using paper-and-pencil tests having slightly better outcomes than fully computer-based classes using only computer assessments.



## 1. Introduction

The use of computers in education and assessment is well established in practice and a widely researched topic. Results accumulate separately by subject and implementation. We examine an important subject to science, technology, engineering, and mathematics (STEM) education. For decades educators and policy makers have sought ways to increase the rate at which the U.S. produces graduates in STEM fields (National Science & Technology Council, 2018). Mathematics—and particularly, calculus—continues to be a persistent gatekeeper to STEM careers. As noted by researchers, failure in an undergraduate mathematics course is one of the greatest predictors of a student's decision to switch to a non-STEM major (Sithole, Chiyaka, McCarthy, Mupinga, Bucklein & Kibirige, 2017). In fact, failure in the calculus sequence often spells the end of a student's efforts to continue in a STEM major (Dibbs, 2019).

Calculus II is usually part of a STEM sequence in the U.S. called mainstream calculus by the Conference Board of the Mathematical Sciences (Blair, Kirkman, & Maxwell, 2018) and the Mathematical Association of America (Bressoud, Mesa, & Rasmussen, 2015). The 2015 Conference Board of the Mathematical Sciences survey estimated that Calculus II had a (nondistance learning) enrollment in fall 2015 of roughly 157,000 students in U.S. higher education (Blair et al., 2018, p. 17, Table S.5). The enrollment in spring may have been even larger since mainstream Calculus I, the feeder course, had a fall 2015 enrollment of 317,000 students (Blair et al., 2018, p. 17, Table S.5). Online homework systems are in widespread use in calculus (Burn & Mesa, 2015, p. 50).

Our focus is on computer-based tests (CBT) versus paper-and-pencil tests (PPT), i.e., the mode of testing. Typical calculus PPT questions may be complex requiring substantial student



work.  Can answer-only CBT evaluate the same knowledge and skill as PPT?  If one accepts the

validity of PPT, then can they be replaced by CBT?  It is a fundamental question if computer

usage is to expand from homework systems to testing.

What is the effect on student learning?  Related current research in other disciplines is

mixed, often suggesting that test mode can make a difference in learning in some disciplines and

for some students (Randal, Sireci, Li, & Kaira, 2012; Prisacari & Danielson, 2017; Vispoel,

Morris, & Clough, 2019).  In calculus, if a class employs a computer homework system and

CBT, then it is fully computer based.  Are there student learning consequences in a fully

computer-based calculus class in which students never submit written work?  Our literature

review did not find any statistical comparison of fully computer-based calculus to one with some

graded pencil-and-paper student work.  We explore computer-based testing as well as fully

computer-based calculus with 324 students.

One reason to use CBT and fully computerized online courses is the advantage of

automated grading over laborious human grading.  Computer systems offer immediate feedback

as well as flexibility for students.  Other reasons may arise.  As we revise this article, we are

under stay-at-home orders due to the COVID-19 crises.   The corresponding author was teaching

PPT calculus classes that were converted to CBT calculus classes out of necessity.

## 2.  Background

### 2.1. Cognitive demand in calculus tasks

We use the framework of White and Mesa's (2014) delineation of cognitive demand,

which is expanded from Tallman and Carlson (2012).  White and Mesa prefer the term

"cognitive orientation" in their analysis of 4,953 calculus tasks partly because "we cannot

observe students' actual interaction with the task" (p. 678).  It is true that we cannot observe



actual interaction, but in examining paper-and-pencil student work, an instructor attempts to glimpse that interaction. White and Mesa's specialized categories are summarized with definitions and examples in their Table 1 (White and Mesa, 2014, p. 680). These task categories are put into a hierarchy of three broad groups: Simple Procedures, Complex Procedures, and Rich Tasks.

Rich tasks are the most complex and make up nearly 50% of calculus test according to White and Mesa's study. They included all tasks coded as Understand, Apply Understanding, Analyze, Evaluate, and Create. However, in practice Rich Tasks are only Understand and Apply Understanding since Analyze, Evaluate, and Create made up less than 0.1% of tasks (White and Mesa, 2014). Whether Understand requires written work is a matter of interpretation. As White and Mesa point out, Tallman and Carlson define Understand as "demonstrate[s] a student has assimilated a concept into an appropriate scheme (Tallman and Carlson 2012, p. 10)" but Krathwohl explicitly requires "determining the meaning of instructional messages, including oral, written, and graphic communication" (Krathwohl 2002, p. 215). Apply Understanding tasks are multistep tasks involving inferences and procedures.

There are two issues we mention concerning grading a complex problem on CBT only-answer criteria in contrast to a PPT implementation that is graded for partial credit. The first is whether obtaining the correct answer is indicates the student mastered the process and understands the concept. Is the answer enough information to render a positive assessment? The second is whether an answer-only CBT question is too high stakes leveling every minor error or misunderstanding to a zero measurement. One can attempt to mitigate these concerns using past experience or best practices. For example, the test should not have answers easily



guessed by trial and error. For the second, inform the student an attempt was wrong, so they may review their scratch work for errors and resubmit.

## 2.2 Computer-Based Tests

Computers have been part of the testing landscape for decades. Their role has evolved from tools for automated scoring of multiple-choice items to be a primary means by which assessments are administered to examinees (Hoogland and Tout, 2018). A persistent question has been how comparable CBT are to their PPT counterparts (Schoenfeldt, 1989). There is substantial research on cognitive workload, studies of internal reliability of various computer-based assessments, and how computer-based assessment fits into developments in learning science compared to paper-and-pencil assessment, e.g., adaptive learning, learner-centered approach (Noyes, Garland, & Robbins, 2004; Prisacari & Danielson, 2017; Emerson & MacKay, 2011; Shute & Rahimi, 2017). There are studies comparing the mode effect on essentially identical assessments. For example, Bennett, Braswell, Oranje, Sandene, Kaplan, and Yan (2008) studied K-12 students taking nearly the same test, question for question. They examined differences in median and variance of scores, and item response theory for difficulty, discrimination, and guessing. In higher education, Clariana and Wallace (2002) studied a college business class, Computer Fundamentals. Students took an identical multiple-choice test by either CBT or PPT. The CBT had a higher median. Hoogland and Tout (2018) observe conflicting pressures and tensions in the development of computer-scored tests. They agree with arguments that prioritizing psychometrics and item response theory in computer-scored tests encourage "assessments of lower-order knowledge and skills" (pp. 680–681). Nevertheless, it is necessary for CBT use in a STEM calculus class that the CBT does assess higher-order knowledge and skills.



**2.2.1. CBT and PPT in calculus.**  In calculus where items may entail proofs, graphic devices and higher order thinking skills, the potential for a testing modality impact could be significant.  Rønning (2017) examined calculus for engineers using Maple TA, a similar system to WebAssign.  His observations were based on student surveys and interviews.  Rønning reports that with PPT students spend time preparing a response for the instructor, thus focusing their energies on the process more than obtaining the correct answer.  In contrast, with CBTs they spend their effort focused on obtaining the correct answer with less attention or focus on the process that was used.  He argues that in upper level mathematics (e.g., calculus) this is particularly problematic where process and proofs have greater importance than in lower level courses.  This is highlighted by the important role partial credit plays in assessment in mathematics.  Partial credit is a means by which instructors can provide guidance to students during the learning process, perhaps rewarding a correct process while not ignoring a calculation error that led to an incorrect answer.  CBTs may attempt to mimic this with answer-until-correct type items, but Rønning claims, this fundamentally changes the nature of the task, reducing it to a "hunting for the correct answer" task with increasingly fewer incorrect options (2017, p. 101).

**C. Student performance in online courses**. In our design, students taking CBT were fully computer based with no pencil-and-paper work to evaluate.  One might expect the mode of testing could affect the raw scores of CBT and PPT.  However, such an affect may be due to the mode of testing affecting student's knowledge and skills rather than reflecting the cognitive workload of tests.  The vehicle for this phenomena would likely be a change in preparation behavior of students (Bennett, 2014).  The course under study uses computer-based homework in all groups and the groups differ only in testing modes.  Therefore, the CBT group is fully computer based while groups using some PPT are not.



The literature on comparisons of computer-based and paper-and-pencil courses in calculus include Hirsch and Weibel (2003), Halcrow and Dunnigan (2012), LaRose (2010), Smolinsky, Olafsson, Marx, and Wang (2019), and Zerr (2007). They all allow that computer-based homework may improve student achievement. However, none included comparisons with fully computer-based courses that exclude all paper-and-pencil graded material. Students were still required to prepare for PPT.

**2.2.2. WebAssign**. Computer homework systems in calculus courses are evenly divided among MyMathLab, WebAssign, and WeBWorK (Kehoe, 2010, p. 755). WebAssign (WebAssign 2019), used for both homework and CBT tests, has copyrighted questions tied to the textbook. Answers may be algebraic, numerical, or multiple choice. WebAssign includes an answer-until-correct feature that allows instructors to set the allowed number of attempts. There is research which suggest WebAssign has a positive effect on faculty and student perceptions of student achievement (Deaveans & Jackson, 2009). In contrast, others suggest little effect on mathematical habits of mind (Nelson & Seminellii, 2016) or metacognitive learning strategies (Khan, 2018). More details of the homework system and the implementation in this study are described in the section *Online homework* of Smolinsky et al. (2019).

**2.2.3. Reliability and validity**. Content validity is nuanced when measuring attributes (e.g., beliefs, attitudes), but less ambiguous when the domain is well defined (DeVellis, 2012, p. 60). What we refer to as a test (either PPT or CBT) consists of series of four tests each covering 3 weeks of class material and a total of 4 chapters of Stewart's textbook (2011). Our tests were functioning class topic tests. DeVellis tells us that an assessment has "content validity when its items are a randomly chosen subset of the universe of appropriate items" (2012, p. 60). Our items were random, but not completely. When the test topic was Taylor series, a question would



be to compute the Taylor series of a function, but the function would be randomized.  Tests reflect a selection of homework from the textbook, which are mainstream.  By 2014 Stewart's Calculus was "used globally by about 70 per cent of students of mathematics, science, and engineering" (Baldassi, 2014, para. 3).  The PPT are mainstream calculus tests graded with partial credit.  Our test questions fall within the structure of White and Mesa's categories (2014) delineating cognitive demand in calculus tests.  It is discussed in the section *Test Items* of Smolinsky et al. (2019, p. 1520). The sample questions in Table 1 of Smolinsky et al. (2019, p. 1520) may be taken as analogous items for our Test $A_p$ and Test $B_p$.  Hence , we believe PPT has content validity to measure knowledge and skill as is covered in mainstream calculus for STEM majors at the level of Stewart's textbook.  This notion of knowledge and skill is our latent variable.

The CBT has the same content questions as PPT, but the grading is different.  The difference in grading may call into question the validity of the CBT. DeVellis's (2012) gives a baseball example, we give another.  Two people may evaluate a batter.  One may look only at hitting accuracy and the second may evaluate the form of the swing by present theory.  Are they evaluating the same thing?  More directly, consider a typical question on both CBT and PPT: Evaluate the integral $\int p^5 \ln(p)\,dp$. For any credit on a CBT, a student must submit an answer algebraically equivalent to the family $C - p^6/36 + (1/6)p^6 \ln(p)$ within the number of allowed submissions.  Logical argument is not explicitly required and if it is measured, then it is indirectly.  In contrast, PPT awards credit for advancing on the logic and these steps must be manifestly exhibited and graded.  The final correct answer on PPT is worth little.  The PPT answer would be a full-page essay/calculation.  An error in one step may set the student on path



to the wrong answer but might be graded near full credit. Hence a PPT and CBT with identical questions may be measuring different latent variables.

We approach the validity of CBT as follows. The PPT and CBT manifestly have the same content. We examine CBT from the standpoint of empirical comparison with the PPT (our gold standard) to establish criterion-related or predictive validity (DeVellis, 2012, p. 61). Construct validity is also established by correlation with PPT—again accepting the validity of PPT (DeVellis, 2012, p. 64).

Reliability has two aspects in this article. The first is the reliability of an instrument. A reliable instrument is an instrument that performs in a consistent manner (DeVellis, 2012, p. 31). Reliability is examined for both CBT and PPT using a standard statistic, Cronbach's alpha (DeVellis, 2012). The second aspect of reliability is the question of whether the CBT and PPT (as two different judges) will give a common assessment and display good inter-rater agreement. We measure the inter-rater reliability of PPT and CCT following model 3 of Shrout and Fleiss (1979) and their statistic ICC(3,2). Model 3 is the model for judges as two different instruments. It is also restated by Gwet as model 10.1 with the comment that "the reliability experiment may use two measuring instruments" as the judges (Gwet, 2004, p. 270).

**2.2.4. Regression and ANCOVA**. CBT and PPT scores may be consistent measures but numerically different. Take an analogous example: Math ACT scores which run from 1 to 36 and the Math SAT which run from 200 to 800. These two raters (instruments) may validly rate the same high school math knowledge and skills, and they may have a strong intraclass correlation indicating that they reliability give a consistent judgement as two different judges. Our PPT and CBT are analogous. A score of 60% may be very good for one instrument and a score of 40% may be the equivalent very good on the other. A strong ICC or Pearson correlation



indicates a linear relationship between scores of the instruments (but not identical numerical or percentile scales). Two instruments might be valid measuring the same underlying variable and might be reliable giving consistent evaluation but have a low correlation. It would mean the conversion of scores is not linear.

For a deeper examination into actual student learning, one must understand how to convert CBT into PPT scores. We compute "the best" linear conversion of scores by the standard of least squares or using linear regression.

In order to measure student learning, one needs a common measure for both students in Section 1 (PPT) and students in Section 3 (CBT) for which to compare using a t-test or ANCOVA (Hinkelmann & Kempthorne, 2007). The common measure is required because a t-test or ANCOVA analyze the difference of means of populations, controlling for other variables, e.g. Gender. The scores from CBT and PPT must be made directly comparable by being put on the same scale using the regression results (from Section 2), otherwise any effect due to Section would be confounded with the different test scales. A t-test or ANCOVA can then be done on the comparable (converted) scores. This method is called method A in this article. Method B involves all the data gathered and several multiple regressions to increase the statistical power.

### 2.3. Research Questions

In our research questions we explore both the consistency of PPT/CBT results and student learning outcomes. Specifically, the two research questions addressed in this study were:

- **Question 1:** Do PPT and CBT measure knowledge and skill in STEM Calculus II (as defined in Section 2.2.3) in a consistent manner? If the answer is yes, then what score y on the PPT indicates the same level of knowledge and skill as a score of x on CBT?



Question 1 is concerned with the reliability and validity of PPT and CBT as well as the inter-rater reliability of PPT and CCT. If the answer is yes, to reliability, validity, and inter-rater reliability, then either test may be used to measure knowledge and skill in Calculus II. The ability to convert CBT scores into equivalent PPT scores will allow us to compare skill and knowledge in a CBT class with the skill and knowledge in a PPT class.

- **Question 2:** How does the knowledge and skill gained by students in a fully computer-based Calculus II class compare to students in a class requiring pencil-and-paper tests.

Students in the CBT section were a fully computer-based class. All homework and tests were computer-based graded on answers only without written work, explanation, or argument ever directly assessed. Students in the PPT section had computer-based homework but were required to complete paper-and-pencil tests requiring work, explanation, or argument. Hence, they also had to study for PPT.

In answering Question 2, we also examined underrepresented minority students. We believe it is advisable to check for effects on minority populations in examining education protocols—even if none is anticipated. We believe it is particularly important to know if a protocol is helpful or harmful in STEM gateway course.

### 3. Method and design

The hosting department offered five sections of its mainstream Calculus II class in spring 2019 that were open enrollment. Students, by their own choice, registered for a section based on their schedule, instructor preference, time preference, and available space. Three sections were used for this experiment. The three sections in this study were listed in the schedule book as team taught by Smolinsky and Olafsson. Smolinsky and Olafsson delivered the lectures, assembled the homework assignments, and authored the tests. Classes met four days a week and



were taught in lecture format with specific days reserved for review. In addition, there was a tutoring center available for students. This is typical of Calculus II as discussed in the study by Smolinsky et al. (2019) and Blair, Kirkman, and Maxwell (2013, p. 17). All classes received the same lectures and used the same textbook (Stewart, 2011). A lecture that was given to one section was also given to all other sections by the same teacher on the same day. Sample tests, sample problems, or handouts were identical for all sections and posted on the course management system, Moodle. Homework was assigned in WebAssign. All three sections were given the same assignments and due dates. The classroom format (except for identical homework and testing format) was used in Smolinsky et al. (2019) and is discussed in more detail there. The difference between treatment of the three sections in this study was in the testing procedures. CBT were given in a monitored and proctored university computer-based testing facility.

### 3.2. Tests

PPT primarily consisted of free answer questions, but sometimes included multiple-choice questions. A multiple-choice item might be included in a two-part question. For example: Does a series converge, and if yes, then find the limit. Determining convergence is a two-choice multiple-choice question. CBT assessments included free answer items and occasionally a multiple-choice item. CBT items had pooled problems of various types and identical problems had randomized parameters since students did not take assessments together.

PPT were hand graded by graduate assistances who award partial credit based on Olafsson's rubric. CBT questions were graded right or wrong by computer. However, because students may get immediate feedback, partial credit was awarded for free answer items: Full credit for a correct answer on the first attempt, 90% on the second, 80% on the third, 70% on the



fourth, 60% on the fifth, 50% on the sixth and last possible attempt.  This method of partial credit was described as "best practice" in WebAssign teacher instructions (WebAssign, 2019).

### 3.3. Sample and Population

The population for this study was primarily undergraduate engineering and science majors at Louisiana State University.  The Colleges of Engineering, Science, and Coast and Environment accounted for 97% of students with declared majors.  While 13.4% had not yet declared a college, as freshmen are not required to declare.

There is a strong preference for randomized controlled trials (Coalition for Evidence-Based Policy, 2003), and a preference for studies in realistic settings (Schoenfeld, 2008).  Actual randomization in actual classes is often not possible.  While the pool of students was not randomly assigned, it was the case that students were unaware that they were participating in a study and students were unaware of testing options in each section.  Students who switched sections (for any reason) after that information became available were eliminated from the study. We have previously referred to this selection process as "virtually randomized."  A discussion of the strengths and weaknesses of virtually randomization are given in the *Treatment and Control Group Selection* section of Smolinsky et al. (2019).  The measures were prospectively arranged, and the only data gathered could have been obtained in nonexperimental sections. This last item was necessary to obtain the university's institutional review board ethics approval without making the students aware of the study.

Participants in the study are students who completed the class and took all exams under the proper settings.  Students who were given special accommodations due to disability, illness, or excused absence on any test were not included, as were students who dropped the class prior to taking the final exam.  The design included controls for Mathematics ACT scores and if the



student received a Pell grant.  If that information was not available, then the student was not included.  The study ultimately included a total of 324 students of which 48 were underrepresented minorities (i.e., African American, Hispanic or Latino, American Indian or Alaska Native, or Pacific Islander), 113 were women, and 86 had Pell grants.

### 3.4. Research Design

The experimental design involved three sections of students.  Each section had four in-semester testing periods that covered material introduced since the previous testing period.  Each of the testing periods required a two-part test covering material that we refer to as part A and part B with each part covering half the material.  The subscript p or c indicates a test respectively is PPT or CBT.  Table 1 outlines the testing design.  Note that $A_p$ was administered to both sections 1 and 2, and $B_c$ was administered to both sections 2 and 3.

**Table 1**

*Testing design*

| Section | Tests | | Final |
|---|---|---|---|
| | part A material | part B material | |
| 1 | $A_p$ | $B_p$ | $F_p$ |
| 2 | $A_c$ | $B_p$ | |
| 3 | $A_c$ | $B_c$ | $F_c$ |

Students' standardized test scores[1] (Math ACT), gender, and Pell grant statuses were investigated as covariates.  Pell grants are the primary source of U.S. government student aid for college students, and awards are based on a computation of the expected family contribution (U.S. Department of Education, 2015). It is a ubiquitous, if imperfect, proxy for socioeconomic status in U.S. higher education (Delisle, 2017).  The Math ACT score (or Math SAT) equivalent was used a pre-test score as a measure of the student's mathematical knowledge at the beginning of the class. Note that Pell and gender are binary variables and Math ACT is treated as continuous.

---

[1] Mathematics ACT score or the Mathematics SAT equivalent score.

**3.4.1. Question 1.**  There are several points to be established in the investigation of Question 1.

a) Find the inter-rater reliability of PPT and CCT using Cronbach's alpha.

b) Measure criterion-related validity (or predictive validity) of CBT using correlation and Cronbach's alpha.

c) Measure the inter-rater reliability of PPT and CCT using ICC(3,2).

d) Compute the conversion of CCT scores to PPT scores.

Items a–c are discussed in Section 2.2.3.  Item d is discussed in 2.2.4 and is not used to answer Question 1.  However, item d is necessary for Question 2 and its computation is naturally done with item b.

**3.4.2. Question 2**.  To address Question 2, we provide two different analyses.  The first uses a straightforward method of analysis using only data that is comparable.  The second method is more complicated giving somewhat sharper p-values, but the cost is limitations on the method.

***3.4.2.1. Method A***.  We compute three pairwise mean comparisons with t-tests, while controlling for Pell grant status, gender, and Math ACT score:

- The score on $B_p$ that is common to both sections 1 and 2.  This allows for a straightforward comparison via a  t-test, adjusted for the covariates .

- The score on $A_c$ that is common to both sections 2 and 3.  This allows for a straightforward comparison via a t-test, adjusted for the covariates.



- A comparison of sections 1 and 3 require additional finesse since there cannot be common assessment, as one section is entirely PPT and the other is entirely CBT. The result of Question 1 part d of section 3.4.1 is a conversion from scores of $A_c$ to $B_p$-equivalent scores, which was established using Section 2 data only. The scores compared in this t-test are $B_p$ scores in Section 1 with the $B_p$-equivalent of the $A_c$ scores in Section 3.

*3.4.2.2. Method B.* Individual test scores are percentages, and scores are averaged whenever test scores are combined. The averaging for this method takes place on the scores after an adjustment by using three multilinear transformations. All scores from tests in Table 1 are converted into $B_p$-equivalent scores using three multilinear transformations arrived at using multiple linear regression with controls for the covariates. We then apply a three-group treatment comparison (PPT, CBT, or Combination PPT/CBT) using an analysis of covariance (ANCOVA), with control variables: Pell grant status, gender, and Math ACT score on the $B_p$-equivalent scores.

We also examined outcomes explicitly on underrepresented minority students. We apply Method B using the same multilinear transformations and a three-group treatment ANCOVA, as it offers the best control on response variance.

## 4. Results

### 4.1. Question 1

We obtain an affirmative answer to Question 1 in part c.

1.a. The statistical evaluation of CBT and PPT is given in Table 2. Results are reported as Cronbach's alpha. The values are well above the so-called minimum acceptable value of 0.70 (Tavakol & Dennick, 2011). The limitation of the comparison is noted in Table 2 under



"Difference":  The material in each chapter was divided into part A and B, so test labeled A and B had somewhat different content while content of  the union of A and B was the same as F, but A$\cup$B and F were administered at different times.

**Table 2**

*Reliability of Assessment Types*

| Protocol | Scores | Difference | Cronbach's alpha | Correlation |
|---|---|---|---|---|
| PPT | $A_p$ and $B_p$ | Content[a] | 0.931 | 0.886 |
| PPT | $(A_p + B_p)/2$ and $F_p$ | Time[b] | 0.943 | 0.897 |
| CBT | $A_c$ and $B_c$ | Content[a] | 0.879 | 0.783 |
| CBT | $(A_c + B_c)/2$ and $F_c$ | Time[b] | 0.871 | 0.789 |

[a]Assessments have different content given at same time.

[b]Assessments have same content given at different times.

1.b.  Cronbach's alpha was 0.886, which finds a respectable value for criterion-related validity above the minimum acceptable value of 0.70.  The alternative measure of correlation was 0.796.

1.c.  ICC(3,2) is equivalent to Cronbach's alpha (Shrout & Fleiss, 1979, p. 422).  It was 0.886 and above the minimum acceptable value of 0.70 (Tavakol & Dennick, 2011).  This gives an affirmative answer to Question 1.

1.d.  With an affirmative answer to Question 1, we may compute how to convert scores. The result of the simple linear regression is shown in Figure 1, with correlation between $A_c$ and $B_p$ being 0.796.  We recognize from the graph that there are possibly issues with low and high outliers, some violation of homoscedasticity, or that a quadratic model may be better.  However, (1) comparisons between the difference of material and in sections 1 and 3 (Table 2) were adequately modeled by a linear relationship meaning the possible deviation would be in the change from PPT to CBT rather than material; (2) we were looking for a change in scale and the



correlation was consistent the results in sections 1 and 3; and (3) we did not want to overfit the data with the introduction of a quadratic.

**Figure 1**

*Simple linear regression of CBT versus PPT*

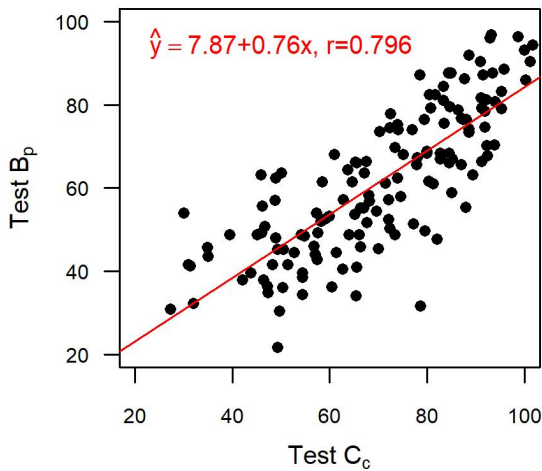

**Each point is an ordered pair of scores for one student in section 2. Correlation (r) and least squares fit for CBT versus PPT.**

We could not expect the correlation in section 2 to be better than the worst of sections 1 and 3. In fact, it could have been significantly worse since the correlation in section 2 was complicated by the reliabilities of both PPT and CBT. The likelihood that actual correlation between $A_c$ and $B_p$ is greater than 0.783 is approximately 0.7, which was estimated using the bootstrap distribution (Zhang & Wang 2017:59). Since we obtain an affirmative answer to Question 1, we use the conversion formula, $y = 0.76 x + 7.87$, to convert the x= $A_c$ score to y = $B_p$- equivalent score.

**4.2. Question 2 Method A**

**4.2.1. B of section 1 versus B of section 2.** Table 3 shows t-test results across sections with different testing protocols, using a model that controls for Pell grant status, math ACT score,



and gender for identical PPT. Section 1 having only PPT, and section 2 having both PPT and CBT. Section 1 did slightly better on average, however the p-value of 0.715 means the outcome is likely random and not statistically significant.

**Table 3**

*t-test for PPT in Sections 1 & 2*

| Effect | Sections | Difference of Means | Standard Error | n | Pr > |t| |
|---|---|---|---|---|---|
| Testing protocol: PPT vs Combination | 1 vs 2 | 1.0927 | 2.9931 | 182 | 0.7155 |

DF = 177

**4.2.2. C of section 2 versus C of section 3.** Table 4 shows t-test results across sections with different testing protocols for identical CBT, using a model that controls for Pell grant status, math ACT score, and gender. Section 2 having both PPT and CBT, and section 3 having only CBT. Section 2 did have a higher mean performance of approximately 2.5 points, but a p-value of 0.213, which is short of even being moderately statistically significant.

**Table 4**

*t-test for CBT in Sections 2 & 3*

| Effect | Sections | Difference of Means | Standard Error | n | Pr > |t| |
|---|---|---|---|---|---|
| Testing protocol: Combination vs CBT | 2 vs 3 | 2.5310 | 2.0298 | 281 | 0.2135 |

DF = 276

**4.2.3. $B_p$ scores from section 1 versus B-equivalent scores of $A_c$ from section 3.** Table 5 shows t-test results for identical CBT across sections with different testing protocols, using a model that controls for Pell grant status, math ACT score, and gender. Section 1 had only PPT, and section 3 had only CBT. The function y = 0.76 *x* +7.87 is the conversion found in investigating



Question 1.  The result was that the PPT section 1 did about 2.77 percentage points better on average when compared to section 3, but the p-value was again non-significant at 0.274.

**Table 5**

*t-test for PPT versus CBT*

| Effect | Sections | Difference of Means | Standard Error | n | Pr > |t| |
|--------|----------|---------------------|----------------|---|---------|
| Testing protocol: PPT vs CBT | 1 vs 3 | 2.6902 | 2.4521 | 182 | 0.2741 |

DF = 180

## 4.3. Question 2 Method B

The results of the ANCOVA with controls for Pell grant status, math ACT score, and gender are given in Table 6.  The comparison of the fully CBT section 3 and the fully PPT section 1 shows a 3.7 point mean preference for PPT in student learning.  The p-value is now moderately significant at 0.084.  There are mean preferences for PPT over both the combined (PPT/CBT) and the combined over the CBT sections, but the non-significant p-values are 0.377 and 0.221, respectively.  The complete ANCOVA is consistent with  PPT improving student learning.  This conclusion is also consistent with the results from method A.  While many of the p-values are somewhat high, there is an overall consistency to the results (in both methods A and B) that favors PPT even though results do not surpass a bright-line rule of p<0.05 (Wasserstein & Lazar, 2016).  The effect size was small.  The standardized effect size (Cohen's-d) between the Section 1 and 3 (the only pair with moderate significance) testing protocols was 1.67. For the testing protocol effect, the 95% confidence intervals for eta-squared and partial eta-squared are (0, 0.036) and (0, 0.039), respectively.



**Table 6**

*ANCOVA*

| Effect | Sections | Difference of Means | Standard Error | Pr > \|t\| |
|---|---|---|---|---|
| Testing protocol: PPT vs Combination | 1 vs 2 | 1.9170 | 2.1658 | 0.3768 |
| Testing protocol: PPT vs CBT | 1 vs 3 | 3.7392 | 2.1557 | 0.0838 |
| Testing protocol: Combination vs CBT | 2 vs 3 | 1.8222 | 1.4851 | 0.2207 |

DF = 324

The results of the ANCOVA for minorities with controls for Pell grant status, math ACT score, and gender are given in Table 7.  There is a preference for CBT section 3 over the combined PPT/CBT section 2 of 8.58 points on average.  This mean difference is substantial in magnitude, while having a moderately significant p-value of 0.063, approaching 0.05.  There is a mean preference for the CBT section 3 over the PPT section 1, but with a p-value greater than 0.9 that is not statistically significant.  The mean of PPT section 1 did better than the combined PPT/CBT section 2, again with a non-compelling p-value of 0.1657.  The results are not consistent with respect to CBT versus PPT.  Further investigation seems to be required with an eye toward the possibility that CBT may be preferred for underrepresented minorities.  A preference for CBT would be surprising as it contradicts the results for the general population. The effect size was small.  The standardized effect size (Cohen's-d) between the Section 2 and 3 (the only pair with moderate significance) testing protocols was 1.98. For the testing protocol effect, the 95% confidence intervals for eta-squared and partial eta-squared are (0, 0.221) and (0, 0.228), respectively.



**Table 7**

*ANCOVA on Underrepresented Minorities*

| Effect | Sections | Difference of Means | Standard Error | Pr > \|t\| |
|---|---|---|---|---|
| Testing protocol: PPT vs Combination | 1 vs 2 | 8.0168 | 5.6824 | 0.1657 |
| Testing protocol: PPT vs CBT | 1 vs 3 | -0.5622 | 4.9399 | 0.9099 |
| Testing protocol: Combination vs CBT | 2 vs 3 | -8.5789 | 4.4885 | 0.0628 |

DF = 42

## 4.4. Effects of Covariates

We find some consistent and interesting covariate effects for Math ACT score, gender, and Pell grant status.

Recall that Tables 3–5 summarized three pairwise mean score comparisons (Method A), comparing PPT vs Combination, Combination vs CBT, and PPT vs CBT, respectively. Although all of these mean differences were found to be non-significant, Table 8 below further provides the p-values associated with the covariate effects. There are three main findings: (a) the Math ACT is highly significant and strongly linearly associated with scores (p-value < 0.0001), (b) The Pell effect is non-significant (p-value > 0.26 in all cases), and (c) The gender effect is either significant or moderately significant, but for only for the results involving the full CBT or section 3, resulting in higher scores for females than males. Specifically, for the PPT vs CBT, females had an average score of about 4.80 higher than males (p-value < 0.053), whereas with Combination vs CBT, females gained 4.30 over males on average (p-value < 0.025).

Tables 6–7 additionally summarized the pairwise comparisons for the testing protocol pairwise average scores comparisons (Method B), comparing PPT vs Combination, Combination vs CBT, and PPT vs CBT, respectively. The companion Table 9 below further presents the p-values associated with the covariate effects. The same main findings are found in Table 9, as



outlined in Table 8: (a) The Math ACT is highly significant and strongly linearly associated with averages (p-value < 0.0001), but less so with the minority analysis (p-value < 0.014) (b), The Pell effect is non-significant (p-value > 0.60 for both the all student and minority analyses), and (c) The gender effect is significant for the all-student analysis (p-value < 0.011), with again females outperforming males on average by 3.72 points. However, the gender effect was not significant for the minority analysis (p-value > 0.28). There were no significant gender by testing protocol interactions found (p-value > 0.80 and 0.42, respectively).

**Table 8**

*Effects of Covariates for Tables 2-4 (Method A)*

| Effect | Sections | Math ACT[*] | Pell[*] | Gender[*] |
|---|---|---|---|---|
| Testing protocol: PPT vs Combination | 1 vs 2 | <0.0001 | 0.5041 | 0.3661 |
| Testing protocol: PPT vs CBT | 1 vs 3 | <0.0001 | 0.2685 | 0.0531 |
| Testing protocol: Combination vs CBT | 2 vs 3 | <0.0001 | 0.9314 | 0.0248 |

[*]$\text{Pr} > |t|$

**Table 9**

*Effects of Covariates for Tables 5-6 (Method B)*

| Subjects | Math ACT[*] | Pell[*] | Gender[*] |
|---|---|---|---|
| All Students  (Model: Table 5) | <.0001 | 0.6018 | 0.0107 |
| Minority Students (Model: Table 6) | 0.0136 | 0.6016 | 0.2812 |

[*]$\text{Pr} > |t|$

## 5.  Discussion

ICC(3,2) or Cronbach's alpha compared the rating of each student over the entire semester and fell into the high good range.  We suspected that CBT and PPT measured different aspects of knowledge and skill or at least directly measured different aspects.  But as David L. Streiner  comments, "α measures not only the homogeneity of the items, but also the



homogeneity of what is being assessed" (2003, p. 102).  It may be that ability to express a logical argument in paper and pencil is indirectly expressed in the ability to think logically and clearly enough to obtain the correct answer to a complicated question.  However, it may also be that the requirements of logical argument in a non-honors STEM calculus class on paper-and-pencil assessments are not particularly strict as—even for mathematics majors—proof and argumentation is usually not taught until after calculus (see the discussion in Smolinsky et al., 2019).  Our results do not indicate that Rønning (2017) warnings on CBT are so widespread as to give a low inter-rater reliability.  His observation that students hunt for the correct answer needs to be taken into account in CBT test construction by avoiding simple numerical answers.

The good inter-rater agreement has some surprising aspects.  For example, the high cognitive demand tasks (in White and Mesa's characterization (2014)) of drawing graphs.  They fall into Krathwohl's Understand (2002), and diverse versions of this task occur throughout the course:  Polar graphs, conic sections, parametric curves in space, quadratic surfaces in space.  In WebAssign, "sketch the graph" is reduced to multiple choice.  In principle, a student may use the process of elimination without being able to actually produce a sketch.  Nevertheless, while the multiple-choice question may be less demanding, it still requires aspects of Understanding such as interpreting, classifying, and comparing in order to eliminate incorrect choices.

Are CBT questions too high stakes?  We expect high-stakes on individual questions to manifest as low reliability.  While the reliability of CBT well above the minimum acceptable value, it was below the PPT value.  Both PPT and CBT scores compared the whole semester, but there was a larger variance for individual complex problems.  The two chapters calculus students historically considered the most difficult are *techniques of integration* and *sequences and series*. For the first test in section 2 that consisted of three complex techniques of integration problems,



the correlation between CBT and PPT was only 0.6 much lower than the overall semester

correlation of 0.796 (Figure 1).

    Four students who got a zero on the CBT portion of test 1 were over 50% on the PPT

portion.  There were no students in the class who received a zero on the PPT portion.  Beyond

the observation that a computer is a merciless grader, is the observation that some students may

do well on partial credit but be unable to correctly complete any multi-step questions.  We did

anecdotally observe (there were no surveys or interviews) the loss of confidence referred to in

Rønning (2017) via some complaints from students who did poorly on the CBT but well on PPT.

We did also hear from students who did poorly on the PPT but well on CBT.  We would not call

it a loss of confidence—they were confident their PPT score should have been graded higher.

There was not a lack of feedback observed by Rønning in computer-based homework.  Students

could go to free tutors, TAs, or instructors to clarify the details of the computer feedback.  The

details of feedback on computer work in this manner still came more quickly than from feedback

via paper and pencil on tests in this study and graded homework in Smolinsky et al. (2019).

    In terms of student outcomes, there were no statistically significant outcomes.

Nevertheless, the results given in Tables 3–6 are all consistent and leaning to paper and pencil.

We were particularly aware that section 3 was fully computer based without students being

required to give paper-and-pencil explanation at any point.  Contrarily, the complexity of the

problems would require scratch work.  To get credit after the first attempted answer would

require students to look over their work for errors, and the complexity of the problems would not

allow guessing as a plausible strategy.  Add to the consistent trend of the results, that the

comparison of the PPT section (1)  and the fully computer section (3) yielded moderately

significant $p<0.084$ in the ANCOVA, we think the results give mild support to PPT.  The point



difference between these two sections was 3.7/100.  Our results are potentially an understatement of the true effects related to test type, as we must approximate $B_p$-equivalents using several regression models leading to an inflation of the experimental error.

In terms of the controls, students with higher math ACT scores did better on tests and women generally did better than men.  There were not significant interactions meaning that math ACT scores and gender did not mean a preference for CBT or PPT.

The results of Table 7 indicate that combination CBT/PPT section did worse than either PPT or CBT sections for minority students.  However, the sample size of 48 was small spread over three sections.  Furthermore, we do not have a theory why the combination CBT/PPT section would be worse than both PPT and CBT individually.  It may be a point for further exploration, but we did attach significance to the results for minority students.

**Limitations of the study**.

Experimenter bias is a limitation whenever the authors conduct a study.  The possibility of experimenter bias was limited by the following procedures:  The outcome measures and covariates were chosen prospectively; the subjective (partial credit) grading was not performed by the authors; the composition of instruments followed the composition for previous non-experimental years prior to conceiving the experiment; the lecture notes and examples were either written for previous non-experimental years prior to conceiving the experiment; or (for some examples) were taken from the textbook.

One reviewer raised the issue of student talking to each other about their testing experiences.  We expected communication to occur and mitigated its effect as best as possible, as outlined below.  In fact, it is an issue whenever flexible scheduling is employed for CBT, but flexible scheduling is often cited as an advantage of CBT.  We made the expectations of the



assessments transparent to all students prior to the test with sample problems.  We would have done so in any (non-experimental) calculus course too.  PPT was simultaneously taken by all students.  While CBT was not simultaneously taken, there were pooled questions and multiple tests when pooling was not technically possible. Hence, we attempted to make communication irrelevant.

The methodology of the ANCOVA (method B) has a limitation in interpreting comparisons with Section 2.  This limitation is because Section 2 is used compute the conversion of scores between CBT to PPT.  Hence the conversion of scores is not independent of the ANCOVA analysis.

## 6.  Conclusion

CBTs are used in many different disciplines and it is likely that their use will continue to grow (Shute & Rahimi, 2017).  They should have advantages for instructors, but our WebAssign implementation of scoring and pooling as well as the availability of questions created time consuming issues for instructors.  Issues that were not relevant for WebAssign's use as a homework system.  When the technical issues are settled, there is the question of PPT versus CBT as assessment tools.  While we found the results favor PPT over CBT, the effect size was small indicative of the little difference in overall performance, across sections.  Furthermore, this was our first implementation, and computer-based assessments may improve with instructor experience and technical development.

### Acknowledgements

The investigators are grateful to the people who assisted in the study. Graduate assistants Phuc A. Nguyen, Bao T. Pham, Federico Salmoiraghi, Zhumagali Shomanov, and Iswarya Sitiraju from the Louisiana State University Department of Mathematics graded paper-and-pencil



exams. The authors are grateful to Eugene Kennedy of the LSU School of Education for pointing them to some of the literature. The manuscript was improved by suggestions from Susan Dunham and the journal's anonymous reviewers.